\documentclass[11pt,reqno]{amsart}
\usepackage{amsmath,amssymb,geometry,color,tikz-cd}
\usepackage[backref,pagebackref,pdftex,hyperindex]{hyperref}
\usepackage[alphabetic]{amsrefs}
\usepackage{amssymb}% see geometry.pdf on how to lay out the page. There's lots.
\usepackage{bbm}
\usepackage{enumitem}
\usepackage{upgreek}

\newtheorem{proposition}{Proposition}[section]
\newtheorem{theorem}[proposition]{Theorem}
\newtheorem{lemma}[proposition]{Lemma}

\theoremstyle{definition}
\newtheorem{remark}[proposition]{Remark}
\newtheorem{definition}[proposition]{Definition}

\title{On wall-crossing for K-stability}
\author{Chuyu Zhou}
\address{\'Ecole Polytechnique F\'ed\'erale de Lausanne (EPFL), MA C3 615, Station 8, 1015 Lausanne, Switzerland}
\email{chuyu.zhou@epfl.ch}

\date{} % delete this line to display the current date

\newcommand{\Fut}{{\rm{Fut}}}

\newcommand{\ord}{{\rm {ord}}}
\newcommand{\tc}{{\rm {tc}}}
\newcommand{\vol}{{\rm {vol}}}

\newcommand{\red}{{\rm {red}}}

\newcommand{\dt}{{\rm {dt}}}

\newcommand{\SL}{{\rm {SL}}}

\newcommand{\PGL}{{\rm {PGL}}}

\newcommand{\Aut}{{\rm {Aut}}}

\newcommand{\bA}{\mathbb{A}}

\newcommand{\bC}{\mathbb{C}}

\newcommand{\bP}{\mathbb{P}}
\newcommand{\bQ}{\mathbb{Q}}

\newcommand{\mD}{\mathcal{D}}
\newcommand{\mE}{\mathcal{E}}
\newcommand{\mF}{\mathcal{F}}
\newcommand{\mG}{\mathcal{G}}

\newcommand{\mL}{\mathcal{L}}
\newcommand{\mM}{\mathcal{M}}

\newcommand{\mO}{\mathcal{O}}

\newcommand{\mW}{\mathcal{W}}
\newcommand{\mX}{\mathcal{X}}
\newcommand{\mY}{\mathcal{Y}}

\begin{document}

\begin{abstract}
In this paper, we explore the wall crossing phenomenon for K-stability, and apply it to explain the wall crossing for K-moduli stacks and K-moduli spaces.

\end{abstract}

\maketitle
\tableofcontents

\section{Introduction}\label{sec: section 1}

By the works \cite{Jiang20, BLX19, Xu20, BX19, ABHLX20, BHLLX21, CP21, XZ20b, LXZ21}, a compactified moduli space is constructed for Fano varieties via the concept of K-stability. Then it is natural to study the relations between different compactifications. For example, the work \cite{ADL19} establishes a wall crossing result for K-moduli of smoothable log Fano varieties, which indicates that K-moduli is likely to play a role of bridge connecting GIT-moduli and KSBA-moduli. One reason for the smoothable condition being posed there is caused by the lack of properness for K-moduli at that time, and they heavily depend on the analytic tools developed in \cite{CDS15a,CDS15b,CDS15c, Tian15, TW20, LWX19} to overcome this difficulty. Along with the establishment of properness by the work \cite{LXZ21} recently, we can clarify the wall crossing phenomenon for K-moduli by a purely algebraic way, and this is the motivation of the current paper.

Let $\pi: \mX\to B$ be a $\bQ$-Gorenstein flat family of $\bQ$-Fano varieties of dimension $d$ such that the anti-canonical divisor $-K_{\mX/B}$ is a relatively ample $\bQ$-line bundle on $\mX$, and $B$ is a normal base. Suppose $\mD$ is an effective $\bQ$-divisor on $\mX$ such that every component of $\mD$ is flat over $B$ and $\mD\sim_\bQ -K_{\mX/B}$ over $B$.
For a rational number $0<c<1$, we denote by $p_c: (\mX, c\mD)\to B$ (the subscript of $p_c$ is used here to distinct the coefficient $c$) and  consider the following subset of the closed points of $B$:
$$B_c:=\{t\in B| \textit{$(\mX_t, c\mD_t)$ is K-semistable}\}. $$
Let us write $B^{ss}\subset B$ to be the subset parametrizing K-semistable fibers of $\pi=p_0$, then $B_0=B^{ss}$, which is an open subset of $B$ (see \cite{BLX19,Xu20}). For the family $p_1: (\mX, \mD)\to B$, let $U\subset B$
be the following subset parametrizing log canonical fibers
$$U:=\{t\in B| \textit{$(\mX_t, \mD_t)$ is log canonical}\}, $$ 
then by inversion of adjunction, $U$ is an open subset of $B$. By interpolation property of K-stability (e.g. \cite[Proposition 2.13]{ADL19}), we see that $B^{ss}\cap U\subset B_c$ for any rational $0<c<1$. Therefore, $B_c$ are mutually birational as $c$ varies for rational $c\in [0,1]$ when $B^{ss}$ and $U$ are both non-empty. Note here that when $c=1$, the pair $(\mX_t,\mD_t)$ being K-semistable just means that it is log canonical (e.g. \cite[Section 9]{BHJ17} or \cite[Theorem 1.5]{Oda13}).
Our goal is to explain the following wall-crossing phenomena for K-stability.
\begin{enumerate}
\item Find finite rational numbers $0=c_0<c_1<c_2<...<c_k<c_{k+1}=1$ such that $B_c$ does not change as $c$ varies in $(c_j,c_{j+1}), j=0,1,...,k$.
\item Figure out the relationships among $B_{c_j-\epsilon}, B_{c_j}$ and $B_{c_j+\epsilon}$ for $0<\epsilon\ll1$ and $1\leq j\leq k$.
\end{enumerate}

This goal is achieved in Theorem \ref{thm: wall-crossing}, and we also apply it to explain the wall crossing phenomenon for K-moduli,  which culminates in Theorem \ref{wall crossing}. 

\begin{theorem}
Let $\pi: \mX\to B$ be a $\bQ$-Gorenstein flat family of $\bQ$-Fano varieties of dim $d$ such that $-K_{\mX/B}$ is a relatively ample $\bQ$-line bundle on $\mX$, and $B$ is a normal base. Suppose $\mD$ is an effective $\bQ$-divisor on $\mX$ such that every component of $\mD$ is flat over $B$ and $\mD\sim_\bQ -K_{\mX/B}$ over $B$.
For a rational number $0\leq c\leq 1$, denote $p_c: (\mX, c\mD)\to B$ and 
$$B_c:=\{t\in B| \textit{$(\mX_t, c\mD_t)$ is K-semistable}\}. $$
Then there exist finite rational numbers $0=c_0<c_1<c_2<...<c_k<c_{k+1}=1$ such that $B_c$ does not change as c varies in $(c_j,c_{j+1})$ for $0\leq j\leq k$.
\end{theorem}

The above theorem is contained in Theorem \ref{thm: wall-crossing}, where we also characterize there the relationships among $B_{c_j-\epsilon}, B_{c_j}$ and $B_{c_j+\epsilon}$ for $0<\epsilon\ll1$ and $1\leq j\leq k$. The strategy  is  to study the following two thresholds and prove their constructible property (see Section \ref{sec: section 3}):
\begin{align}
l(t):=&\inf \{0\leq c\leq 1| \textit{\rm the pair $(\mX_t,c\mD_t)$ is K-semistable}\},\\
u(t):=&\sup \{0\leq c\leq 1| \textit{\rm the pair $(\mX_t,c\mD_t)$ is K-semistable}\}.
\end{align}
Note that such kinds of thresholds have been defined for family of smoothable Fano varieties in the works \cite[Definition 7.4]{LWX19} and \cite[Definition 3.14]{ADL19}, and some constructible properties are also proved therein to get openness property for K-stability. In this paper, we will study their constructibility by a different method, which fits to the recently developed algebraic K-stability theory.

In order to apply it to the wall crossing for K-moduli, we first set the notation.
Fix a rational number $r>0$, we consider a set $\mF$ of pairs satisfying that $(X,D)\in \mF$ if and only if it is  of the following form: 
\begin{enumerate}
\item $X$ is a $\bQ$-Fano variety of dimension $d$ and volume $v$ (i.e. $(-K_X)^d=v$),
\item $D\sim_\bQ -rK_X$ is a Weil divisor on $X$ such that $(X, \frac{c}{r}D)$ is K-semistable for some rational $c\in [0,1)$.
\end{enumerate}
In Theorem \ref{bounded lemma}, we will see that $\mF$ lies in a log bounded family. Then for sufficiently large $m$, we define the following subset of $\mF$, denoted by $\mF_{m,\bP^{N_m}}^{\chi,\tilde{\chi}}$  and abbreviated by $\mF_m$, such that $(X,D)\in \mF_m$ if and only if the following conditions are satisfied:
\begin{enumerate}
\item $(X,D)$ can be embedded into $\bP^{N_m}\times \bP^{N_m}$ such that 
$$\mO_X(1):=\mO_{\bP^{N_m}}(1)|_X\sim -mK_X\quad and \quad H^0(\bP^{N_m}, \mO_{\bP^{N_m}}(1))\cong H^0(X,-mK_X),$$
\item $X$ (resp. $D$) admits Hilbert polynomial $\chi$ (resp. $\tilde{\chi}$)
with respect to $\mO_X(1)$ (resp. $\mO_D(1):=\mO_{\bP^{N_m}}(1)|_D$).
\end{enumerate}
Write $P:=(\chi,\tilde{\chi})$. We define the subset $\mF_{m,c}\subset \mF_m$ for a rational number $c\in [0,1)$ as follows:$$\mF_{m,c}:=\{(X,D)\in \mF_m| \textit{\rm the pair $(X,\frac{c}{r}D)$ is K-semistable}\}. $$
By the proof of \cite[Theorem 2.21]{XZ20b}, the set $\mF_m$ admits an algebraic structure, and $\mF_{m,c}\subset \mF_m$ is an open sub-scheme by \cite{BLX19,Xu20}, then we write $\mF_m^{\red}$ and $\mF_{m,c}^{\red}$
to be the corresponding reduced schemes. Denote by 
$\mM^K_{P,r,c}:=[\mF_{m,c}^{\red}/\PGL(N_m+1)]$, which is a reduced Artin stack of finite type representing the moduli pseudo-functor over reduced bases (by the works \cite{Jiang20, BLX19, Xu20}). Here one can also define the moduli functor over any base due to \cite{Kollar19}, see Remark \ref{moduli functor}. We denote by $M^K_{P,r,c}$ the descent of $\mM^K_{P,r,c}$, which is a projective separated good moduli space  parametrizing K-polystable elements, whose existence is due to the works \cite{ABHLX20, BHLLX21, XZ20b, LXZ21}. Then we have the following result.

\begin{theorem}
There exist finite rational numbers $0=c_0<c_1<c_2<...<c_k<c_{k+1}=1$ such that 
\begin{enumerate}
\item the schemes $\mF_{m,c}^{\red}$ do not change as $c$ varies in $(c_j,c_{j+1})$ for $0\leq j\leq k$, and we have open immersions
$$\mF^{\red}_{m,c_j-\epsilon}\hookrightarrow \mF^{\red}_{m,c_j}\hookleftarrow \mF^{\red}_{m,c_j+\epsilon}$$
for rational $0<\epsilon\ll1$ and $1\leq j\leq k$. 
\item the Artin stacks $\mM^K_{P,r,c}$  do not change as $c$ varies in $(c_j,c_{j+1})$ for $0\leq j\leq k$, and we have the following open immersions 
$$\mM^K_{P,r,c_j-\epsilon}\hookrightarrow \mM^K_{P,r,c_j}\hookleftarrow  \mM^K_{P,r,c_j+\epsilon}$$
for rational $0<\epsilon\ll1$ and $1\leq j\leq k$. 
\item  the good moduli spaces $M^K_{P,r,c}$ do not change as $c$ varies in $(c_j,c_{j+1})$ for $0\leq j\leq k$, and we have the following projective morphisms among good moduli spaces
$$M^K_{P,r,c_j-\epsilon}\stackrel{\phi^-_j}{\rightarrow}M^K_{P,r,c_j}\stackrel{\phi^+_j}{\leftarrow}M^K_{P,r,c_j+\epsilon}$$
for rational $0<\epsilon\ll 1$ and $1\leq j\leq k$, which are induced by the following diagram:
\begin{center}
	\begin{tikzcd}[column sep = 2em, row sep = 2em]
	 \mM_{P,r, c_j-\epsilon}^K \arrow[d,"",swap] \arrow[rr,""]&& \mM_{P,r, c_j}^K\arrow[d,"",swap] &&\mM_{P,r, c_j+\epsilon}^K\arrow[d,""]\arrow[ll,""]\\
	 M_{P,r, c_j-\epsilon}^K\arrow[rr,"\phi_j^-"]&& M_{P,r, c_j}^K&&M_{P,r, c_j+\epsilon}^K\arrow[ll,"\phi_j^+", swap].
	 	 	 	\end{tikzcd}
\end{center}
.
\end{enumerate}
\end{theorem}

\noindent
\subsection*{Acknowledgement}
   The author would like to thank Yuchen Liu and  Ziquan Zhuang for helpful discussions and suggestions, and Zsolt Patakfalvi for beneficial comments. The author is supported by grant European Research Council (ERC-804334).

\section{Preliminaries}

In this section, we review the concept of K-stability of Fano varieties via test configurations and valuative criterion. We work over complex number field $\bC$ through the paper. We say that $(X,\Delta)$ is a log pair if $X$ is a projective normal variety and $\Delta$ is an effective $\bQ$-divisor on $X$ such that $K_X+\Delta$ is $\bQ$-Cartier. We say that a log pair $(X,\Delta)$ is log Fano if it admits klt singularities and $-K_X-\Delta$ is ample. If a log Fano pair $(X,\Delta)$ has vanishing boundary, i.e. $\Delta=0$, we just say that $X$ is a $\bQ$-Fano variety. For various concepts of singularities in birational geometry, such as klt, plt, lc, etc., we refer to \cite{KM98, Kollar13}.

\subsection{K-stability via test configurations}

The definition of K-stability via test configurations originates to \cite{Tian97} and was later algebraically reformulated by \cite{Don02}.

\begin{definition}
Let $(X,\Delta)$ be a log pair of dimension $d$ and $L$ an ample $\bQ$-line bundle on $X$. A test configuration $\pi: (\mX,\Delta_\tc;\mL)\to \bA^1$ is a degenerating family over $\bA^1$ consisting of the following data:
\begin{enumerate}
\item $\pi: \mX\to \bA^1$ is a projective flat morphism from a normal variety $\mX$, $\Delta_\tc$ is an effective $\bQ$-divisor on $\mX$, and $\mL$ is a relatively ample $\bQ$-line bundle on $\mX$,
\item the family $\pi$ admits a $\bC^*$-action which lifts the natural $\bC^*$-action on $\bA^1$ such that $(\mX,\Delta_\tc; \mL)\times_{\bA^1}\bC^*$ is $\bC^*$-equivariantly isomorphic to $(X, \Delta; L)\times_{\bA^1}\bC^*$.
\end{enumerate}
\end{definition}

We denote $(\bar{\mX}, \bar{\Delta}_\tc;\bar{\mL})\to \bP^1$ to be the natural compactification of the original test configuration, which is obtained by glueing $(\mX, \Delta_\tc;\mL)$ and $(X,\Delta;L)\times (\bP^1\setminus 0)$ along their common open subset $(X, \Delta;L)\times \bC^*$.

Suppose $(X,\Delta)$ is a log Fano pair of dimension $d$ and $L=-K_X-\Delta$. Let  $(\mX,\Delta_\tc; \mL)$ be a test configuration such that $\mL=-K_{\mX/\bA^1}-\Delta_\tc$. We call it a special (resp. weakly special) test configuration if  $(\mX, \mX_0+\Delta_{\tc})$ admits plt (resp. lc) singularities. When we come to weakly special test configurations, we always assume that the central fibers are integral.

\begin{definition}\label{tc}
The generalized Futaki invariant of a test configuration $(\mX,\Delta_\tc;\mL)$ is defined as follows:
$$\Fut(\mX,\Delta_\tc;\mL):=\frac{(K_{\bar{\mX}/\bP^1}+\bar{\Delta}_\tc).\bar{\mL}^d}{L^d} -\frac{d}{d+1}\frac{\{(K_X+\Delta)L^{d-1}\}\bar{\mL}^{d+1}}{(L^d)^2}.$$
In particular, if $(X,\Delta)$ is a log Fano pair and $L=-K_X-\Delta$,  we have
$$\Fut(\mX,\Delta_\tc;\mL)=\frac{d\bar{\mL}^{d+1}}{(d+1)(-K_X-\Delta)^d} +\frac{\bar{\mL}^d.(K_{\bar{\mX}/\bP^1}+\bar{\Delta}_\tc)}{(-K_X-\Delta)^d}.$$
\end{definition}

\begin{definition}
Let $(X,\Delta)$ be a log Fano pair of dimension $d$ and $L=-K_X-\Delta$. We say that $(X,\Delta)$ is K-semistable if the generalized Futaki invariants are non-negative for all test configurations.
We say that $(X,\Delta)$ is K-polystable if it is K-semistable and any test configuration with K-semistable central fiber is induced by a one parameter subgroup of $\Aut(X,\Delta)$, which is also called a test configuration of product type. 
\end{definition}

\subsection{Valuative criterion}

Let $(X,\Delta)$ be a log Fano pair of dimension $d$, we say that $E$ is a prime divisor over $X$ if there is a proper birational morphism $f: Y\to X$ between normal varieties such that $E$ appears as a prime divisor on $Y$. We say that $E$ is a special (resp. weakly special) divisor over $X$ if $\ord_E$ is induced by some non-trivial special (resp. weakly special) test configuration of $(X,\Delta;L:=-K_X-\Delta)$. In other words, if $(\mX,\Delta_\tc;\mL)$ is a non-trivial special (resp. weakly special) test configuration of $(X,\Delta;L)$, by \cite[Lemma 4.1]{BHJ17}, the central fiber induces a divisorial valuation on $K(X)$, i.e. $\ord_{\ord_{\mX_0}}|_{K(X)}=c\cdot\ord_E$ for some rational number $c>0$ and some prime divisor $E$ over $X$, and such a prime divisor $E$ is called a special (resp. weakly special) divisor over $X$.

\begin{definition}
Let $E$ be a prime divisor over $X$, we have the following invariants:
$$A_{X,\Delta}(E)=\ord_E(K_Y-f^*(K_X+\Delta))+1 ,$$
$$S_{X,\Delta}(E)=\frac{1}{\vol(-K_X-\Delta)}\int_0^\infty\vol(-f^*(K_X+\Delta)-tE)\dt.$$
Then beta invariant is defined to be 
$$\beta_{X,\Delta}(E)=A_{X,\Delta}(E)-S_{X,\Delta}(E);$$ 
and delta invariant is defined to be 
$$\delta(X,\Delta)=\inf_F \frac{A_{X,\Delta}(F)}{S_{X,\Delta}(F)},$$
where $F$ runs through all prime divisors over $X$.
\end{definition}

\begin{remark}\label{rk:supplement remark}
Notation as in the above definition, we give a few remarks here,
\begin{enumerate}
\item The invariant $S_{X,\Delta}(E)$ can also be defined as the limit of $S_m(E):=\max_{D_m}\ord_E(D_m)$, where $D_m$ runs through all $m$-basis type divisors. Note that, if we denote $N_m:=\dim H^0(X,-m(K_X+\Delta))$, then the divisor of the form 
$$\frac{\sum_{j=1}^{N_m}{\rm div}(s_j=0)}{mN_m}\sim_\bQ -(K_X+\Delta)$$ is called an $m$-basis type divisor, where $\{s_j\}_{j=1}^{N_m}$ is a complete basis of  the vector space $H^0(X,-m(K_X+\Delta))$.
\item Define $T_{X,\Delta}(E):=\sup_{0\leq D\sim_\bQ -K_X-\Delta}\ord_E(D)$, then 
we have $\frac{d+1}{d}\leq \frac{T_{X,\Delta}(E)}{S_{X,\Delta}(E)}\leq d+1$, see \cite{BJ20}.
\item For a weakly special test configuration $(\mX,\Delta_\tc;\mL)$ of $(X,\Delta;L)$, denote by $v_{\mX_0}:=\ord_{\mX_0}|_{K(X)}=c\cdot \ord_E$ for some rational number $c>0$ and prime divisor $E$ over $X$, then we have $\Fut(\mX,\Delta_\tc;\mL)=c\cdot\beta_{X,\Delta}(E)$, see \cite[Theorem 6.13]{Fuj19}. 
\item If $\delta(X,\Delta)<\frac{d+1}{d}$, then the $\inf$ in the definition of $\delta(X,\Delta)$ can be replaced by $\min$ due to the recent work \cite{LXZ21}.
\end{enumerate}
\end{remark}

The following theorem is now a cornerstone in the recently developed algebraic K-stability theory.

\begin{theorem}{\rm (\cite{Fuj19,Li17, BJ20,FO18})}
Let $(X,\Delta)$ be a log Fano pair of dimension $d$, then we have the following equivalences:
\begin{enumerate}
\item $(X,\Delta)$ is K-semistable if and only if $\beta_{X,\Delta}(E))\geq 0$ for any prime divisor $E$ over $X$ (or equivalently for any special divisor $E$ over $X$).
\item $(X,\Delta)$ is K-semistable if and only if $\delta(X,\Delta)\geq 1$ for any prime divisor $E$ over $X$ (or equivalently for any special divisor $E$ over $X$).

\end{enumerate}
\end{theorem}

\subsection{Fano type varieties}

Let $f: Y\to Z$ be a projective morphism such that $Y$ is a normal variety. We say that $Y$ is Fano type over $Z$ if there exists an effective $\bQ$-divisor $B$ on $Y$ such that $(Y, B)$ has klt singularities and $-K_Y-B$ is ample over $Z$. When $Z$ is a point, we just say that $Y$ is of Fano type or $Y$ is a Fano type variety.

A Fano type variety is naturally a Mori dream space (e.g. \cite{BCHM10}). This means that for a given projective morphism $Y\to Z$ such that $Y$ is Fano type over $Z$, one can run $D$-MMP over $Z$ for 
any $\bQ$-Cartier divisor $D$ to get a minimal model or a Mori fiber space.

\subsection{Log Calabi-Yau pairs}

We say that a log pair $(Y,G)$ is log Calabi-Yau if it admits log canonical singularities and $K_Y+G\sim_\bQ 0$. In this subsection, we only focus on log Calabi-Yau pairs of the form $(X,D)$, where $X$ is a $\bQ$-Fano variety of dimension $d$ and $D\sim_\bQ -K_X$. 
Given a log Calabi-Yau pair $(X,D)$, suppose $(\mX, \mD;\mL)\to \bA^1$ is a test configuration of $(X,D;L:=-K_X)$, by Definition \ref{tc}, we have the following formula:
$$\Fut(\mX,\mD;\mL)=\frac{(K_{\bar{\mX}/\bP^1}+\bar{\mD})\bar{\mL}^d}{L^d}. $$
For $c\in [0,1]$, we also have the following computation. 

\begin{lemma}\label{interpolation computation}
Notation as above, we have
$$\Fut(\mX, c\mD;\mL)=c\cdot\Fut(\mX,\mD;\mL) +(1-c)\cdot \Fut(\mX;\mL).$$
\end{lemma}

\begin{proof}
The equality is clear for $c=0,1$. We assume $c\in (0,1)$. By Definition \ref{tc}, 
\begin{align*}
\Fut(\mX,c\mD;\mL)=&\ \frac{(K_{\bar{\mX}/\bP^1}+c\bar{\mD}).\bar{\mL}^d}{L^d} -\frac{d}{d+1}\frac{\{(K_X+cD)L^{d-1}\}\bar{\mL}^{d+1}}{(L^d)^2}\\
=&\ \frac{(K_{\bar{\mX}/\bP^1}+c\bar{\mD}).\bar{\mL}^d}{L^d}+\frac{d(1-c)}{d+1}\frac{\bar{\mL}^{d+1}}{L^d}\\
=&\ (1-c)\cdot \{ \frac{K_{\bar{\mX}/\bP^1}.\bar{\mL}^d}{L^d} + \frac{d}{d+1}\frac{\bar{\mL}^{d+1}}{L^d} \}+c\cdot \frac{(K_{\bar{\mX}/\bP^1}+\bar{\mD})\bar{\mL}^d}{L^d}\\
=&\  (1-c)\cdot \Fut(\mX;\mL) +c\cdot\Fut(\mX,\mD;\mL).
\end{align*}
\end{proof}

We will need the the following lemma later.
\begin{lemma}\label{lemma on cy}
Notation as above, we have the following two statements:
\begin{enumerate}
\item Suppose $E$ is an lc place of $(X,D)$, then one can produce a test configuration $(\mX,\mD; \mL=-K_{\mX/\bA^1})$ of $(X,D;-K_X)$  such that  $(\mX, \mX_0+\mD)$ is log canonical with $\mX_0$ being integral and 
$$\ord_{\mX_0}|_{K(X)}=\ord_E \quad and \quad \Fut(\mX,\mD;\mL)=0.$$
\item Suppose $(\mX,\mD; \mL=-K_{\mX/\bA^1})$ is a test configuration of $(X,D;-K_X)$ such that $(\mX,\mX_0+ \mD)$ is log canonical with $\mX_0$ being integral, then $\ord_{\mX_0}|_{K(X)}$ is an lc place of $(X,D)$.
\end{enumerate}
\end{lemma}

\begin{proof}
We first prove (1). 
By \cite[Corollary 1.4.3]{BCHM10}, there exists an extraction morphism $f: Y\to X$ which only extracts $E$. If $E$ is a prime divisor on $X$, we take $f={\rm id}$. We just assume that $E$ is exceptional over $X$, as the case when $E$ is a prime divisor on $X$ can be addressed similarly.  We have the following:
$$K_Y+f_*^{-1}D+E=f^*(K_X+D). $$
First observe that $Y$ is of Fano type, so $E$ is a dreamy divisor (e.g. \cite[Corollary 1.3.2]{BCHM10} and \cite[Definition 1.3]{Fuj19}).
Put the morphism $f$ into the constant family $F: Y_{\bA^1}\to X_{\bA^1}$, then we have
$$K_{Y_{\bA^1}}+F_*^{-1}D_{\bA^1}+E_{\bA^1}+Y_0=F^*(K_{X_{\bA^1}}+D_{\bA^1}+X_0) ,$$
where $X_0$ (resp. $Y_0$) is the central fiber of the family $X_{\bA^1}\to \bA^1$ (resp. $Y_{\bA^1}\to \bA^1$). 
Let $v$ be a divisorial valuation over $X_{\bA^1}$ with weights $(1,1)$ along divisors $E_{\bA^1}$ and $Y_0$ (e.g. \cite[2.3.2]{BLX19}), then 
$$A_{X_{\bA^1}, D_{\bA^1}+X_0}(v)=A_{Y_{\bA^1}, F_*^{-1}D_{\bA^1}+E_{\bA^1}+Y_0}(v)=0. $$
By \cite[Corollary 1.4.3]{BCHM10}, one can extract $v$ to be a divisor on a birational model $g: \mY\to X_{\bA^1}$ as follows:
$$K_{\mY}+g_*^{-1}D_{\bA^1}+g_*^{-1}X_0+\mE=g^*(K_{X_{\bA^1}}+D_{\bA^1}+X_0),$$
where $v=\ord_\mE$.
The pair $(\mY, g_*^{-1}D_{\bA^1}+g_*^{-1}X_0+\mE)$ is  log canonical and being log Calabi-Yau over $\bA^1$ (i.e. $K_{\mY}+g_*^{-1}D_{\bA^1}+g_*^{-1}X_0+\mE\sim_\bQ 0$ over $\bA^1$), and $\mY$ is of Fano type over $\bA^1$.
Then we could run $g_*^{-1}X_0$-MMP over $\bA^1$ to get a minimal model $\mY\dashrightarrow \mY'/\bA^1$. Denote $D_{\mY'}$ (resp. $X_0'$ and  
$\mE'$) to be the pushforward of $g_*^{-1}D_{\bA^1}$ (resp. $g_*^{-1}X_0$ and $\mE$). By Zariski lemma (e.g. \cite[Section 2.4]{LX14}), $g_*^{-1}X_0$ is contracted, that is, $X_0'=0$. Thus the following two pairs
$$(\mY, g_*^{-1}D_{\bA^1}+g_*^{-1}X_0+\mE) \dashrightarrow (\mY', D_{\mY'}+\mE'),$$
are both log canonical and being log Calabi-Yau over $\bA^1$, and $\mY'$ is of Fano type over $\bA^1$. 
We focus on the pair $(\mY', D_{\mY'}+\mE')\to \bA^1$. Run $-K_{\mY'}(\sim_{\bQ,\bA^1} D_{\mY'}+\mE')$-MMP over $\bA^1$ to get an ample model $\mY'\dashrightarrow \mX$. Let  $\mD$ (resp. $\mE''$) be the pushforward of $D_{\mY'}$ (resp. $\mE'$), then we see that $(\mX,\mD) \to \bA^1$ is a test configuration of $(X,D)$, and the central fiber is exactly $\mE''$. By our construction of $\mE$ we know that the restriction of $\ord_{\mE''}$ to $K(X)$ is exactly $\ord_E$.  

Let $p:\mW\to (\bar{\mX},\bar{\mD})$ and $q: \mW\to (X_{\bP^1},D_{\bP^1})$ be a the common log resolution, by above construction we see that $p^*(K_{\bar{\mX}/\bP^1}+\bar{\mD})=q^*(K_{X_{\bP^1}/\bP^1}+D_{\bP^1})\sim_\bQ 0$, which deduces the vanishing of $\Fut(\mX,\mD;\mL)$.

Now we turn to (2). Denote by $\ord_{\mX_0}|_{K(X)}=c\cdot \ord_E$ for some rational $c>0$ and prime divisor $E$ over $X$. Then we see that $(\mX,a\mD)$ is a weakly special test configuration of $(X, aD)$ for any rational $a\in [0,1)$. By \cite[Appendix]{BLX19}, for any rational $0<\epsilon\ll 1$, one can find a $D_\epsilon\sim_\bQ -K_X$ such that $E$ is an lc place of $(X,(1-\epsilon)D+\epsilon D_\epsilon)$. Suppose $E$ is not an lc place of $(X,D)$, then $A_{X,D}(E)>t$ for some $t>0$. On the other hand, 
$$A_{X,(1-\epsilon)D+\epsilon D_\epsilon}(E)=A_{X, (1-\epsilon)D}(E)-\epsilon \cdot\ord_E(D_\epsilon)\geq t-\epsilon \cdot T_X(E).$$
For sufficiently small $\epsilon$ we have $A_{X,(1-\epsilon)D+\epsilon D_\epsilon}(E)>0$, a contradiction.
\end{proof}

\section{Constructibility of delta invariants for family of log Fano pairs}\label{sec; section 2}

Let $(X,\Delta)$ be a log Fano pair of dimension $d$ and $D\sim_\bQ -K_X-\Delta$ an effective $\bQ$-divisor on $X$. Given a real number $0<c<1$, we define 
$$\delta(X,\Delta+cD):=\inf_E \frac{A_{X,\Delta+cD}(E)}{(1-c)S_{X,\Delta}(E)},$$ 
where $E$ runs through all prime divisors over $X$. Note here that $\delta(X,\Delta+cD)$ is exactly the delta invariant of the log Fano pair $(X,\Delta+cD)$ when $c$ is rational. For convenience, we also say that $(X, \Delta+cD)$ is a log Fano pair if it admits klt singularities, though $c\in (0,1)$ is irrational.

\begin{definition}
For a given real number $0<c<1$, we say the pair $(X,\Delta+cD)$ is K-semistable if it is klt and $\delta(X,\Delta+cD)\geq 1$.
\end{definition}

Note that this definition is now standard if $c$ is rational due to the works \cite{FO18,BJ20}. However, if $c$ is irrational, it may not be standard to talk about the concept of K-stability for the pair in the literature, as one should start with test configurations and prove the valuative criterion as in \cite{Li17,Fuj19}. Here we do not need the equivalence between the two definitions of K-stability, respectively via test configurations and via valuative criterion, for real coefficients, and we just define K-stability in a way as the above definition. Define 
$$\tilde{\delta}(X,\Delta+cD):=\min\{\delta(X,\Delta+cD), 1\},$$ 
we have the following generalization of openness property for K-stability in \cite{BLX19,Xu20}.

\begin{theorem}{\rm (\cite{BLX19,Xu20})}\label{real constructibility}
Let $\pi: \mX\to B$ be a $\bQ$-Gorenstein flat family of $\bQ$-Fano varieties of dim $d$ such that $-K_{\mX/B}$ is a relatively ample $\bQ$-line bundle on $\mX$, and $B$ is a normal base. Suppose $\mD$ is an effective $\bQ$-divisor on $\mX$ such that every component is flat over $B$ and $\mD\sim_\bQ -K_{\mX/B}$ over $B$. For any given real number $0<c<1$, the set $\{\tilde{\delta}(\mX_t,c\mD_t)| t\in B\}$ is finite.
\end{theorem}

Note that the result is a direct consequence of \cite{BLX19,Xu20} for rational $c$. For irrational $c$, it is enough to show that $\delta(X,cD)$ can be approximated by lc places of bounded complements of $(X,-K_X)$ if it is not greater than $1$ (see the lemma below), and the rest follows from \cite[Proposition 4.3]{BLX19}.
 \begin{lemma}
 Let $(X,\Delta)$ be a log Fano pair of dimension $d$ and $D\sim_\bQ -K_X-\Delta$ an effective $\bQ$-divisor on $X$. Fix an irrational number $0<c<1$ such that the log Fano pair $(X,\Delta+cD)$ satisfies that $\delta(X,\Delta+cD)\leq 1$. Then there exists an integer $N$ depending only on $d$ and the coefficients of $\Delta$ such that 
 $$\delta(X,cD)=\inf_E \frac{A_{X, \Delta+cD}(E)}{(1-c)S_{X,\Delta}(E)}, $$
 where $E$ runs through lc places of $N$-complements of $(X,\Delta)$.
 \end{lemma}
 
 \begin{proof}
We define 
$$\delta_m:=\delta_m(X,\Delta+cD):=\inf_E \frac{A_{X,\Delta+cD}(E)}{\max_{D_m}(1-c)\ord_E(D_m)},$$ 
 where $E$ runs through all prime divisors over $X$ and $D_m$ runs through all $m$-basis type divisors of $(X,\Delta)$. This infimum is in fact a minimum, since $m$-basis divisors form a bounded family and log canonical thresholds take only finitely many values on this family.
 It is clear that $\lim_m \delta_m(X,\Delta+cD)=\delta(X,\Delta+cD)$ (see Remark \ref{rk:supplement remark} (1)).
 
 We first assume $\delta(X,\Delta+cD)<1$, then $\delta_m(X, \Delta+cD)<1$ for sufficiently large $m$. In this case one can find an $m$-basis type divisor $D_m$ such that $(X,\Delta+cD+ \delta_m(1-c)D_m)$ is strictly log canonical and admits an lc place $E_m$ over $X$. We choose a general $H\sim_\bQ -K_X-\Delta$ such that $(X,\Delta+cD+\delta_m(1-c)D_m+(1-c-\delta_m(1-c))H)$ is an lc log Calabi-Yau pair. By \cite[Corollary 1.4.3]{BCHM10}, there is an extraction $g_m: Y_m\to X$ which only extracts $E_m$, i.e.
 \begin{align*}
 &K_{Y_m}+\tilde{\Delta}+c\tilde{D}+\delta_m(1-c)\tilde{D}_m+(1-c-\delta_m(1-c))\tilde{H}+E_m \\
 =&g_m^*(K_X+\Delta+cD+\delta_m(1-c)D_m+(1-c-\delta_m(1-c))H),
 \end{align*}
 where $\tilde{\bullet}$ is the birational transformation of $\bullet$ and $Y_m$ is of Fano type. We run $-(K_{Y_m}+\tilde{\Delta}+E_m)$-MMP to get a mimimal model $f_m: Y_m\dashrightarrow Y_m'$ such that $-(K_{Y_m}+{f_m}_*\tilde{\Delta}+{f_m}_*E_m)$ is nef and $Y_m'$ is also of Fano type. By \cite[Theorem 1.7]{Birkar19}, there exists a positive integer $N$ depending only on $d$ and the coefficients of $\Delta$ such that $E_m$ is an lc place of some $N$-complement of $(X,\Delta)$. Therefore we have
 $$\delta(X,\Delta+cD)=\lim_m \frac{A_{X,\Delta+cD}(E_m)}{(1-c)S_{X,\Delta}(E_m)},$$
 where $E_m$ are a sequence of lc places of some $N$-complements of $(X,\Delta)$.
 
 Next we turn to the case $\delta(X,\Delta+cD)=1$. By the same proof as \cite[Theorem 1.1]{ZZ19}, for any rational $0<\epsilon_i<1$, one can find a $D_i\sim_\bQ -K_X-\Delta$ such that $\delta(X,\Delta+cD+\epsilon_iD_i)<1$. We assume $\{\epsilon_i\}$ is a decreasing sequence of rational numbers tending to $0$. For each $\epsilon_i$, by the proof of the previous case, there exists a sequence of prime divisors $\{E_{i,m}\}_m$ which are lc places of some $N$-complements of $(X,\Delta)$ such that
 $$\delta(X,\Delta+cD+\epsilon_i D_i)=\lim_m \frac{A_{X,\Delta+(c+\epsilon_i)D}(E_{i,m})}{(1-c-\epsilon_i)S_{X,\Delta}(E_{i,m})},$$ 
and $N$ depends only on $d$ and the coefficients of $\Delta$. Thus we obtain the following 
 $$\delta(X,\Delta+cD)=\inf_{i,m}\frac{A_{X,\Delta+cD}(E_{i,m})}{(1-c)S_{X,\Delta}(E_{i,m})}. $$
 The proof is finished.
 \end{proof}
 
\begin{proof}[Proof of Theorem \ref{real constructibility}]
The proof is just a combination of the above lemma and \cite[Proposition 4.3]{BLX19}. 
\end{proof} 
 
 We also need the following lemma in the next section. 
 \begin{lemma}\label{vanishing futaki}
 Let $(X,\Delta)$ be a log Fano pair and $D\sim_\bQ -K_X-\Delta$ an effective $\bQ$-divisor on $X$. Suppose $0<c<1$ is a given real number such that $(X,\Delta+cD)$ is a K-semistable log Fano pair, and $E$ is a prime divisor over $X$ such that it induces a special test configuration $(\mX, \Delta_{\tc}+c\mD)$ of $(X,\Delta+cD)$. If $(\mX_0, \Delta_{\tc,0}+c\mD_0)$ is K-semistable, then 
 $$A_{X,\Delta+cD}(E)-(1-c)S_{X,\Delta}(E) =0.$$
 \end{lemma}
 
 Note that it is well known if $c$ is rational. 
 \begin{proof}
As $(X,\Delta+cD)$ is K-semistable,  we have $A_{X,\Delta+cD}(E)-(1-c)S_{X,\Delta}(E)\geq 0$. 
 By the same proof as \cite[Theorem 5.1]{Fuj19}, the generalized Futaki invariant of the test configuration $(\mX,\Delta_\tc+c\mD)$ is equal to 
 $A_{X,\Delta+cD}(E)-(1-c)S_{X,\Delta}(E)$. Suppose it is strictly positive, then there exists a product type test configuration of $(\mX_0,\Delta_{\tc,0}+c\mD_0)$ with negative generalized Futaki invariant\footnote{The $\bC^*$-action on $(\mX_0, \Delta_{\tc,0}+c\mD_0)$ induces a product type test configuration of $(\mX_0, \Delta_{\tc,0}+c\mD_0; -(K_{\mX_0}+\Delta_{\tc,0}+c\mD_0))$. By \cite{Wang12, Odaka13b}, the generalized Futaki invariant of the product type test configuration  is the total weight of the $\bC^*$-action on $-(K_{\mX_0}+\Delta_{\tc,0}+c\mD_0)$. If we compose this $\bC^*$-action with $\bC^*\to \bC^*, t\mapsto t^{-1}$, then we get another product type test configuration of $(\mX_0, \Delta_{\tc,0}+c\mD_0; -(K_{\mX_0}+\Delta_{\tc,0}+c\mD_0))$ with opposite total weight.}. This will produce a prime divisor $F$ over $\mX_0$ such that 
 $$A_{\mX_0,\Delta_{\tc,0}+c\mD}(F)-(1-c)S_{\mX,\Delta_{\tc,0}}(F)<0,$$
 which is a contradiction to the K-semistability of $(\mX_0,\Delta_{\tc,0}+c\mD_0)$.
 \end{proof}

\section{Two thresholds $l(t)$ and $u(t)$}\label{sec: section 3}

Our strategy for the goal in the first section is  to study the following two thresholds:
\begin{align}
l(t):=&\inf \{0\leq c\leq 1| \textit{\rm the pair $(\mX_t,c\mD_t)$ is K-semistable}\},\\
u(t):=&\sup \{0\leq c\leq 1| \textit{\rm the pair $(\mX_t,c\mD_t)$ is K-semistable}\}.
\end{align}
In this section, we will study the constructibility of these two thresholds by an approach different from \cite[Definition 7.4]{LWX19} and \cite[Definition 3.14]{ADL19}, which fits to the recently developed algebraic K-stability theory. It is clear that these two functions are defined on the base $B$, where $t\in B$ varies through closed points. It may appear for some $t\in B$ such that the pair $(\mX_t,c\mD_t)$ is not K-semistable for any $c\in [0,1]$, and in this case we just define $l(t)=u(t)=-\infty$. 
We first have the following lemma on semi-continuity of $l(t)$ and $u(t)$.

\begin{lemma}
The function $l(t)$ is upper semi-continuous and $u(t)$ is lower semi-continuous. Therefore, the set $\{l(t)|t\in B\}$ (resp. $\{u(t)|t\in B\}$) satisfies ACC (resp. DCC) condition.
\end{lemma}

\begin{proof}
To show $l(t)$ is upper semi-continuous, it suffices to show that the set 
$$l_a:=\{t\in B| 0\leq l(t)<a\}$$ 
is an open subset for any $0<a<1$. Suppose $t_0\in l_a$, then the pair $(\mX_{t_0},a_0\mD_{t_0})$ is K-semistable for some rational $a_0<a$ by Lemma \ref{lem:rational lemma}(4). By openness property for K-semistable locus (see \cite{BLX19,Xu20}), we know there exists an open subset $U_{t_0}\subset B$ such that the pair $(\mX_t,a_0\mD_{t})$ is K-semistable for any $t\in U_{t_0}$. Thus $l(t)\leq a_0<a$ for any $t\in U_{t_0}$, which concludes the upper semi-continuity of $l(t)$.
The same analysis applies to the lower semi-continuity of $u(t)$. The ACC (resp. DCC) property follows immediately. 
\end{proof}

The ACC (resp. DCC) property of $l(t)$ (resp. $u(t)$) is not enough for our purpose, and we will prove that they even satisfy constructible property. Before that, we will first list the following properties of the two functions.

\begin{lemma}\label{lem:rational lemma}
Notation as above, we have the following properties for $u(t)$ and $l(t)$:
\begin{enumerate}
\item  the pair $(\mX_t,c\mD_t)$ is K-semistable for every rational $c\in [l(t),u(t)]$,
\item  $l(t)=0$ and $u(t)=1$ for $t\in B^{ss}\cap U$ (see the notation in the first section),
\item The $\rm inf$ (resp. $\rm sup$) in the definition of $l(t)$ (resp. $u(t)$) can be replaced by $\rm min$ (resp. $\rm max$),
\item $l(t)$ and $u(t)$ are both rational for any $t\in B$ if they are not $-\infty$,
\item $B_c\subset B$ is an open subset for any rational $0\leq c\leq 1$,
\item $\{t\in B| l(t)=u(t)=-\infty\}$ is a closed subset of $B$.
\end{enumerate}
\end{lemma}

\begin{proof}
For (1), it is a consequence of the combination of (3) with the interpolation property of K-stability (\cite[Proposition 2.13]{ADL19}).  For (2), it is a consequence of interpolation property of K-stability (\cite[Proposition 2.13]{ADL19}). (3) and (4) are corollaries of the next lemma, and (5) is a corollary of openness property of K-semistability by \cite{BLX19,Xu20}.
For (6), it is clear that the subset $\{t\in B| l(t)=u(t)=-\infty\}$ is the complement in $B$ of the following subset:
$$\{t\in B| \textit{$(\mX_t,c\mD_t)$ is K-semistable for some rational $c\in [0,1]$}\},$$
which is an open subset of $B$ by \cite{BLX19,Xu20}, thus (6) is concluded. 
\end{proof}

\begin{lemma}\label{rationality}
Let $X$ be a $\bQ$-Fano variety and $D\sim_\bQ -K_X$ an effective $\bQ$-divisor on $X$. Suppose $0<c<1$ is a real number, then we have the following results:
\begin{enumerate}
\item  If $(X,aD)$ is K-semistable for any $a\in [c-\epsilon,c]$ (resp. $a\in [c,c+\epsilon]$) and some $0<\epsilon \ll1$, but K-unstable for any $a\in (c,1)$ (resp. $a\in (0,c)$), then the number $c$ must be rational.
\item If  $(X,aD)$ is K-semistable for any $a\in [c-\epsilon, c)$ or $a\in (c,c+\epsilon]$ for some $0<\epsilon \ll1$, then $(X,cD)$ is also K-semistable.
\end{enumerate}
\end{lemma}

\begin{proof}
For (1), we only assume $(X,aD)$ is K-semistable for any $a\in [c-\epsilon,c]$ and some $0<\epsilon \ll1$ but K-unstable for any $a\in (c,1)$, as the other case can be obtained by the same way. We choose a sequence of decreasing rational numbers $c<c_i<1$ tending to $c$. By the assumption we see that $(X,c_iD)$ is a K-unstable log Fano pair for each $i$. By \cite{BLZ19,LXZ21}, there is a prime divisor $E_i$ over $X$ such that 
$$A_{X,c_iD}(E_i)-\delta(X,c_iD)S_{X,c_iD}(E_i)=0, $$
and $E_i$ induces a special test configuration of $(X,c_iD)$, denoted by $(\mX_i, c_i\mD_i)$, such that 
$$\delta(\mX_{i,0}, c_i\mD_{i,0})=\delta(X,c_iD).$$
We first show that the set $\{(\mX_{i,0}, \mD_{i,0})\}$ lies in a log bounded family. It is clear  that $\vol(-K_{\mX_{i,0}})=\vol(-K_X)$ and 
$$\frac{A_{\mX_{i,0}}(E)}{S_{\mX_{i,0}}(E)}\geq  (1-c_i)\cdot\frac{A_{\mX_{i,0}, c_i\mD_{i,0}}(E)}{(1-c_i)S_{\mX_{i,0}}(E)}$$
for every prime divisor $E$ over $\mX_{i,0}$. In other words, 
$$\delta(\mX_{i,0})\geq (1-c_i)\delta(\mX_{i,0}, c_i\mD_{i,0})=(1-c_i)\delta(X,c_iD).$$
By the following estimate on $\delta(X,c_iD)$,
\begin{align*}
\delta(X,c_iD)&\ =\ \inf_F \frac{A_{X,c_iD}(F)}{S_{X,c_iD}(F)}=\inf_F \frac{A_{X,cD}(F)-(c_i-c)\ord_F(D)}{\frac{1-c_i}{1-c}(1-c)S_X(F)}\\
&\ \geq\ \inf_F \frac{A_{X,cD}(F)-(c_i-c)T_X(F)}{\frac{1-c_i}{1-c}(1-c)S_X(F)}\\
&\ \geq\  \frac{1-c}{1-c_i}- \frac{(d+1)(c_i-c)}{1-c_i},
\end{align*}
where $F$ runs through all prime divisors over $X$, we see that $\delta(\mX_{i,0})$ admits a positive lower bound which does not depend on $i$ since 
$$\delta(\mX_{i,0})\geq 1-c-(d+1)(c_i-c). $$
Thus the set $\{\mX_{i,0}\}$ lies in a bounded family by \cite{Jiang20}, which implies that the set $\{(\mX_{i,0},\mD_{i,0})\}$ lies in a log bounded family as $\mD_{i,0}\sim_\bQ-K_{\mX_{i,0}}$ and the Weil index of $\mD_{i,0}$ is upper bounded. We now consider the set $\{(\mX_{i,0}, c\mD_{i,0})\}$, which lies in a log bounded family. We have the following estimate on $\delta(\mX_{i,0},c\mD_{i,0})$:
$$\frac{A_{\mX_{i,0},c\mD_{i,0}}(F)}{(1-c)S_{\mX_{i,0}}(F)}= \frac{A_{\mX_{i,0},c_i\mD_{i,0}}(F)+(c_i-c)\ord_F(\mD_{i,0})}{\frac{1-c}{1-c_i}S_{\mX_{i,0}, c_i\mD_{i,0}}(F)}\geq \frac{1-c_i}{1-c}\frac{A_{\mX_{i,0},c_i\mD_{i,0}}(F)}{S_{\mX_{i,0}, c_i\mD_{i,0}}(F)},$$
where $F$ runs through all prime divisors over $\mX_{i,0}$.
Take the infimum on both sides, we have
$$\delta(\mX_{i,0}, c\mD_{i,0})\geq \frac{1-c_i}{1-c}\delta(\mX_{i,0}, c_i\mD_{i,0}). $$
Combine the estimate of $\delta(X,c_iD)=\delta(\mX_{i,0},c_i\mD_{i,0})$ above, we see that
$$\delta(\mX_{i,0}, c\mD_{i,0})\geq 1-  \frac{(d+1)(c_i-c)}{1-c}.$$
Since $(\mX_{i,0}, c\mD_{i,0})$ is log bounded and from Theorem \ref{real constructibility}, we see that $\delta(\mX_{i,0}, c\mD_{i,0})\geq 1$ for $i\gg1$. By Lemma \ref{vanishing futaki}, this implies 
$$A_{X,cD}(E_i)-(1-c)S_X(E_i)=A_X(E_i)-S_X(E_i)-c(\ord_{E_i}(D)-S_X(E_i))=0.$$
By the proof of \cite[Theorem 4,7]{BLZ19} we know that $S_X(E_i)$ is a rational number. Combine the inequality
$$A_{X,c_iD}(E_i)-(1-c_i)S_X(E_i)=A_X(E_i)-S_X(E_i)-c_i(\ord_{E_i}(D)-S_X(E_i))<0,$$
we see that 
$$\ord_{E_i}(D)-S_X(E_i)\ne 0. $$
Thus 
$$c=\frac{A_X(E_i)-S_X(E_i)}{\ord_{E_i}(D)-S_X(E_i)} $$
is a rational number.

For (2), we only assume $(X,aD)$ is K-semistable for any $a\in [c-\epsilon,c)$ for some $0<\epsilon \ll1$, as the other case can be done by the same way. Suppose $(X,cD)$ is K-unstable, then there exists a prime divisor $E$ over $X$ such that 
$$A_{X,cD}(E)-(1-c)S_{X}(E)<0.$$ 
Replace $c$ by another rational $0<c^-<c$ such that $0<c-c^-\ll 1$, one still has 
$$A_{X,c^-D}(E)-(1-c^-)S_{X}(E)=A_{X,c^-D}(E)-S_{X,c^-D}(E)<0,$$ 
which is a contradiction since $(X,aD)$ is K-semistable for any rational $a\in [c-\epsilon,c)$ for some $0<\epsilon \ll1$.
\end{proof}

We are ready to prove the following constructible property.

\begin{theorem}
Both $\{l(t)|t\in B\}$ and $\{u(t)|t\in B\}$ are finite rational sets.
\end{theorem}

\begin{proof}
We first show that $l(t)$ satisfies DCC property. Suppose not, one can find a strictly decreasing sequence $\{0\leq a_i<1\}_{i=1}^{\infty}$ which tends to a limit $a\in [0,1)$, such that there is a corresponding sequence of closed points on $B$, denoted by $\{t_i\}_{i=1}^\infty$, satisfying the following conditions:
\begin{enumerate}
\item The pair $(\mX_{t_i},a_i\mD_{t_i})$ is K-semistable for every $i$,
\item The pair $(\mX_{t_i},c\mD_{t_i})$ is not K-semistable for each rational $0\leq  c<a_i$ and every $i$.
\end{enumerate}
We consider the family $p_a: (\mX,a\mD)\to B$. By the finiteness of the set  (see Theorem \ref{real constructibility})
$$\{\min\{1,\delta(\mX_t,a\mD_t)\}| t\in B\},$$
we see that there exists a positive  number $0<\eta_a<1$ such that, the pair $(\mX_t,a\mD_t)$ is K-semistable if and only if $\delta(\mX_t,a\mD_t)> \eta_a$. Since $(\mX_{t_i},a_i\mD_{t_i})$ is K-semistable, by valuative criterion, 
$$\frac{A_{\mX_{t_i},a_i\mD_{t_i}}(E)}{S_{\mX_{t_i},a_i\mD_{t_i}}(E)} \geq 1$$
for any prime divisor $E$ over $\mX_{t_i}$. Note that $$S_{\mX_{t_i},a\mD_{t_i}}(E)=\frac{1-a}{1-a_i}S_{\mX_{t_i},a_i\mD_{t_i}}(E),$$ 
thus we have 
$$\frac{A_{\mX_{t_i},a\mD_{t_i}}(E)}{S_{\mX_{t_i},a\mD_{t_i}}(E)}\geq  \frac{A_{\mX_{t_i},a_i\mD_{t_i}}(E)}{\frac{1-a}{1-a_i}S_{\mX_{t_i},a_i\mD_{t_i}}(E)}\geq \frac{1-a_i}{1-a} $$
for any prime divisor $E$ over $\mX_{t_i}$. As $a_i$ strictly decreases to $a$, we have $\frac{1-a_i}{1-a}>\eta_a$ for sufficiently large $i$. Thus $(\mX_{t_i},a\mD_{t_i})$ is K-semistable for sufficiently large $i$. This is a contradiction since $(\mX_{t_i},c\mD_{t_i})$ is not K-semistable for any rational $0\leq c<a_i$ and any $i$. The contradiction implies that $l(t)$ satisfies DCC condition, which derives that $\{l(t)| t\in B\}$ is indeed a finite set.

Next we show $u(t)$ satisfies ACC property. Suppose not, one can find a strictly increasing sequence $\{0<a_i< 1\}_{i=1}^{\infty}$ which tends to a limit $a\in (0,1]$, such that there is a corresponding sequence of closed points on $B$, denoted by $\{t_i\}_{i=1}^\infty$, satisfying the following conditions:
\begin{enumerate}
\item The pair $(\mX_{t_i},a_i\mD_{t_i})$ is K-semistable for every $i$,
\item The pair $(\mX_{t_i},c\mD_{t_i})$ is not K-semistable for any rational $a_i<c\leq 1$ and every $i$.
\end{enumerate}
As $(\mX_{t_i},a_i\mD_{t_i})$ is klt for every $i$, by ACC of log canonical thresholds (see \cite{HMX14}), we may assume $(\mX_{t_i},a\mD_{t_i})$ is log canonical for every $i$ after subtracting those $a_i$ that are not sufficiently close to $a$. 

We first assume $a=1$, then there is a positive rational number $0<\xi<1$ such that $(\mX_{t_i},(1-\xi)\mD_{t_i})$ is log canonical if and only if $(\mX_{t_i},\mD_{t_i})$ is log canonical for every $i$. Since $(\mX_{t_i},a_i\mD_{t_i})$ is K-semistable, the pair $(\mX_{t_i},\mD_{t_i})$
is log canonical for sufficiently large $i$, thus being K-semistable by \cite[Theorem 1.5]{Oda13}. By interpolation property of K-stability (\cite[Proposition 2.13]{ADL19}), we see 
$(\mX_{t_i},c\mD_{t_i})$ is K-semistable for any rational $a_i< c\leq 1$, which is a contradiction. 

We next assume $0<a<1$ and consider the family $p_a: (\mX,a\mD)\to B$. 
As we have seen, the set
$$\{\min\{1,\delta(\mX_t,a\mD_t)\}| t\in B\}$$
is finite by Theorem \ref{real constructibility}, thus there exists a positive rational number $0<\xi_a<1$ such that, the pair $(\mX_t,a\mD_{t})$ is K-semistable if and only if $\delta(\mX_t,a\mD_t)> \xi_a$. Since $(\mX_{t_i},a_i\mD_{t_i})$ is K-semistable, by valuative criterion, 
$$\frac{A_{\mX_{t_i},a_i\mD_{t_i}}(E)}{S_{\mX_{t_i},a_i\mD_{t_i}}(E)} \geq 1$$
for any prime divisor $E$ over $\mX_{t_i}$. Again by the relation $$S_{\mX_{t_i},a\mD_{t_i}}(E)=\frac{1-a}{1-a_i}S_{\mX_{t_i},a_i\mD_{t_i}}(E),$$  
we have 
\begin{align*}
\frac{A_{\mX_{t_i},a\mD_{t_i}}(E)}{S_{\mX_{t_i},a\mD_{t_i}}(E)}=\ &\ \frac{A_{\mX_{t_i},a_i\mD_{t_i}}(E)-(a-a_i)\ord_E(\mD_{t_i})}{\frac{1-a}{1-a_i}S_{\mX_{t_i},a_i\mD_{t_i}}(E)} \\
\geq\ &\  \frac{1-a_i}{1-a}-\frac{a-a_i}{1-a}\frac{T_{\mX_{t_i}}(E)}{S_{\mX_{t_i}}(E)}\\
\geq\ &\ \frac{1-a_i}{1-a}-\frac{(a-a_i)(d+1)}{1-a}
\geq \xi_a
\end{align*}
for sufficiently large $i$. Thus the pair $(\mX_{t_i},c\mD_{t_i})$ is K-semistable for each rational $a_i<c\leq a$, which is a contradiction and implies that $\{u(t)|t\in B\}$ is a finite set. The proof is finished.

\end{proof}

As $\{l(t)|t\in B\}$ and $\{u(t)|t\in B\}$ are both finite rational sets, we assume $l(t)$ takes values in $\{0\leq l_i\leq 1| i=0,1,...,r,r+1\}$ with $0=l_0<l_1<l_2<...<l_r<l_{r+1}=1$, and $u(t)$ takes values in $\{0\leq u_j\leq 1| j=0,1,...,s,s+1\}$ with $0=u_0<u_1<u_2<...<u_s<u_{s+1}=1$. We put these rational numbers together and reorder them by identifying the repeating numbers as the following sequence:
$$0=c_0<c_1<c_2<...<c_k<c_{k+1}=1 .$$
We will see in the following theorem that the above sequence of $c_j$ is exactly what we need in Section \ref{sec: section 1}.

\begin{theorem}\label{thm: wall-crossing}
Let $\pi: \mX\to B$ be a $\bQ$-Gorenstein flat family of $\bQ$-Fano varieties of dim $d$ such that $-K_{\mX/B}$ is a relatively ample $\bQ$-line bundle on $\mX$, and $B$ is a normal base. Suppose $\mD$ is an effective $\bQ$-divisor on $\mX$ such that every component of $\mD$ is flat over $B$ and $\mD\sim_\bQ -K_{\mX/B}$ over $B$.
For a rational number $0\leq c\leq 1$, denote by $p_c: (\mX, c\mD)\to B$ and 
$$B_c:=\{t\in B| \textit{$(\mX_t, c\mD_t)$ is K-semistable}\}. $$
Then there exist finite rational numbers $0=c_0<c_1<c_2<...<c_k<c_{k+1}=1$ such that $B_c$ does not change as c varies in $(c_j,c_{j+1})$ for $0\leq j\leq k$.
Moreover, we give the following description of the relationships among $B_{c_j-\epsilon}, B_{c_j}$ and $B_{c_j+\epsilon}$ for rational $0<\epsilon\ll1$ and $0\leq j\leq k$:
\begin{enumerate}
\item $B_0=B^{ss}$ and $B_\epsilon=B^{ss}\setminus \{t\in B| u(t)\leq 0\}$,
\item if $c_j=l_p$ for some $1\leq p\leq r$ and $c_j\ne u_q$ for any $1\leq q\leq s$, we have 
$$B_{c_j-\epsilon}\subset B_{c_j}=B_{c_j+\epsilon}\quad and \quad B_{c_j}\setminus B_{c_j-\epsilon}=\{t\in B| l(t)=c_j\},$$
\item if $c_j=u_q$ for some $1\leq q\leq s$ and $c_j\ne l_p$ for any $1\leq p\leq r$, we have 
$$B_{c_j-\epsilon}= B_{c_j}\supset B_{c_j+\epsilon}\quad and \quad B_{c_j}\setminus B_{c_j+\epsilon}=\{t\in B| u(t)=c_j\},$$
\item if $c_j=l_p=u_q$ for some $1\leq p\leq r$ and $1\leq q\leq s$, we have 
$$B_{c_j-\epsilon}\subset B_{c_j}\supset B_{c_j+\epsilon},$$ 
$$B_{c_j}\setminus B_{c_j-\epsilon}=\{t\in B| l(t)=c_j\}\quad and \quad  B_{c_j}\setminus B_{c_j+\epsilon}=\{t\in B| u(t)=c_j\}.$$
\end{enumerate}
\end{theorem}

\begin{proof}
We first show that $B_c$ does not change as c varies in $(c_j,c_{j+1})$ for $0\leq j\leq k$. Suppose $t\in B_c$, where $c\in (c_j, c_{j+1})$, it suffices to show that $t\in B_a$ for any $a\in (c_j, c_{j+1})$. By our choice of $c_k, k=0,1,..., k+1$, we have $(c_j, c_{j+1})\subset [l(t), u(t)]$, which means that $(\mX_t, a\mD_t)$ is K-semistable for any $a\in (c_j, c_{j+1})$. Thus $t\in B_a$ for any $a\in (c_j, c_{j+1})$. 

Next we show the four statements listed in the theorem. (1) is clear. For (2), suppose $t\in B_{c_j}$ for some $1\leq j\leq k$. Since $c_j\ne u_q$ for any $1\leq q\leq s$, we see that $t\in B_{c_j+\epsilon}$ for $0<\epsilon\ll 1$, i.e. $B_{c_j}\subset B_{c_j+\epsilon}$. Now we suppose $t\in B_{c_j+\epsilon}$ for $0<\epsilon\ll 1$. If $l(t)\ne c_j$, then $t\in B_{c_j}$ obviously; if $l(t)=c_j$, by Lemma \ref{rationality} (2), we still have that $(\mX_t, c_j\mD_t)$ is K-semisatble, i.e. $t\in B_{c_j}$. Thus $B_{c_j}=B_{c_j+\epsilon}$. Suppose $t\in B_{c_j-\epsilon}$ for $0<\epsilon\ll 1$, then we see that $t\in B_{c_j}$ still by Lemma \ref{rationality} (2), thus $B_{c_j-\epsilon}\subset B_{c_j}$ for $0<\epsilon \ll 1$. Suppose $t\in B_{c_j}$, if $l(t)\ne c_j$, then $t\in B_{c_j-\epsilon}$; 
if $l(t)=c_j$, then by the definition of $l(t)$ we see that $t$ is not contained in $B_{c_j-\epsilon}$, thus $B_{c_j}\setminus B_{c_j-\epsilon}=\{t\in B| l(t)=c_j\}$. By now, (2) is proved. The analysis of (3) and (4) are totally the same as that of (2), we just omit the proof.
\end{proof}

We end this section by the following lemma, which will be used in the next section.

\begin{lemma}\label{detailed stability}
Notation as in Theorem \ref{thm: wall-crossing}, we have the following results:
\begin{enumerate}
\item if $c_j=l_p$ for some $1\leq p\leq r$ and $c_j\ne u_q$ for any $1\leq q\leq s$, then the set $\{(\mX_t,c_j\mD_t)|l(t)=c_j\}$ contains no K-polystable element,
\item if $c_j=u_q$ for some $1\leq q\leq s$ and $c_j\ne l_p$ for any $1\leq p\leq r$, then the set $\{(\mX_t,c_j\mD_t)|u(t)=c_j\}$ contains no K-polystable element,
\item for $1\leq j\leq k$ and rational $0<\epsilon\ll 1$, it cannot happen that $(\mX_t, (c_j+\epsilon) \mD_t)$ (resp. $(\mX_t, (c_j-\epsilon)\mD_t)$) is strictly K-semistable (i.e. K-semistable but not K-polystable) but $(\mX_t, c_j\mD_t)$ is K-polystable.
\end{enumerate}
\end{lemma}

\begin{proof}
We only prove (1) and (3), as (2) can be proved by the same way as (1). 
For (1), suppose $t\in B$ satisfies that $l(t)=c_j=l_p$ for some $1\leq p\leq r$ and $c_j\ne u_q$ for any $1\leq q\leq s$, and $(\mX_t,c_j\mD_t)$ is K-polystable. We first see that $(\mX_t, (c_j+\epsilon)\mD_t)$ is K-polystable for all rational $0<\epsilon\ll1$, and $(\mX_t, (c_j-\epsilon)\mD_t)$ is K-unstable for all rational $0<\epsilon\ll1$. By \cite{BLZ19, LXZ21}, there is a special test configuration of $(\mX_t, (c_j-\epsilon)\mD_t)$, denoted by $(\mY_\epsilon, (c_j-\epsilon)\mD_\epsilon)$, such that 
$$\delta(\mX_t,(c_j-\epsilon)\mD_t)=\delta(\mY_{\epsilon, 0}, (c_j-\epsilon)\mD_{\epsilon,0})<1.$$
Thus the generalized Futaki invariant of this test configuration is negative, which implies that it is not of product type\footnote{If the test configuration is of product type, then the generalized Futaki invariant of $(\mY_\epsilon, c_j\mD_\epsilon)$ and that of $(\mY_\epsilon, (c_j+\epsilon)\mD_\epsilon)$ both vanish by the K-polystability of $(\mX_t, c_j\mD_t)$ and $(\mX_t, (c_j+\epsilon)\mD_t)$. By linearity, the generalized Futaki invariant of $(\mY_\epsilon, (c_j-\epsilon)\mD_\epsilon)$ also vanishes, which is a contradiction to its negativity.}. By the same method as the proof of Lemma \ref{rationality} (namely by showing the log boundedness of the set $\{(\mY_{\epsilon,0},c_j\mD_{\epsilon, 0})\}_{0<\epsilon\ll 1}$ and bounding their delta invariants, and then apply Theorem \ref{real constructibility}), we know that $(\mY_{\epsilon, 0}, c_j\mD_{\epsilon,0})$ is K-semistable when $\epsilon$ is sufficiently small, then one sees that the test configuration is of product type, since the K-polystable pair $(\mX_t,c_j\mD_t)$ degenerates to a K-semistable pair $(\mY_{\epsilon,0},c_j\mD_{\epsilon,0})$, contradiction.

For (3), we only assume $(\mX_t, (c_j+\epsilon)\mD_t)$ is strictly K-semistable for rational $0<\epsilon\ll1$ but $(\mX_t, c_j\mD_t)$ is K-polystable, as the other case can be proved similarly. For any special test configuration of $(\mX_t, (c_j+\epsilon)\mD_t)$, denoted by $(\mW_{\epsilon}, (c_j+\epsilon)\mD_\epsilon)$, such that the central fiber $(\mW_{\epsilon,0}, (c_j+\epsilon)\mD_{\epsilon,0})$ is K-semistable, by the same method as the proof of Lemma \ref{rationality}, we know $(\mW_{\epsilon,0}, c_j\mD_{\epsilon,0})$ is also K-semistable for sufficiently small $\epsilon$. Thus the test configuration is of product type, since $(\mX_t,c_j\mD_t)$ is K-polystable. This implies that $(\mX_t, (c_j+\epsilon)\mD_t)$ is K-polystable, contradiction.
\end{proof}

\section{Application to K-moduli}

In this section we study the wall crossing phenomenon for K-moduli. Through the section, we fix a rational number $r>0$. 
We start with the following lemma.
\begin{lemma}\label{lct gap}
Suppose $r\geq 1$. There exists a positive rational number $0<t_0<1$ depending only on $d$ and $r$ such that
if $X$ is a $\bQ$-Fano variety of dimension $d$ and $D\sim_\bQ -rK_X$ is a Weil divisor satisfying that $(X,\frac{1-t_0}{r}D)$ is log canonical, then the pair $(X,\frac{1}{r}D)$ is also log canonical.
\end{lemma}

\begin{proof}
The existence of $t_0$ is just an application of ACC of log canonical thresholds by \cite{HMX14}.
\end{proof}

We also need the following lemma.
\begin{lemma}\label{almost cy}
Notation as in the previous lemma, there exists a rational number $0<t_0'<1$ depending only on $d$ and $r$, such that the pair $(X,\frac{c}{r}D)$ is K-semistable for some rational $c\in [1-t_0',1)$ if and only if $(X,\frac{c}{r}D)$ is K-semistable for any rational $c\in [1-t'_0,1]$.
\end{lemma}
\begin{proof}
It follows from the next lemma by taking $t_0'=\frac{\epsilon_0}{2}$.\end{proof}

\begin{lemma}
Let $(Y,G)$ be an lc log Calabi-Yau pair where $Y$ is a $\bQ$-Fano variety of dimension $d$ and $G\sim_\bQ-K_X$ is an effective $\bQ$-divisor on $Y$. Then there exists a rational $0<\epsilon_0< 1$ depending only on $d$ and the coefficients of $G$ such that the following statements are equivalent:
\begin{enumerate}
\item the pair $(Y,cG)$ for some rational $c\in (1-\epsilon_0,1)$ is K-semistable,
\item the pair $(Y,cG)$ for any rational $c\in (1-\epsilon_0,1]$ is K-semistable,
\item $\beta_Y(E)=A_Y(E)-S_Y(E)\geq 0$ for every lc place $E$ of $(Y,G)$.
\end{enumerate}
\end{lemma}

\begin{proof}
The direction $(2)\Rightarrow (1)$ is clear. We first prove $(1)\Rightarrow (3)$.
Let $E$ be an lc place of $(Y,G)$, then by Lemma \ref{lemma on cy} one can construct a test configuration $(\mY,\mG;-K_{\mY/\bA^1})$ of $(Y,G;-K_Y)$ such that $\Fut(\mY,\mG;-K_{\mY/\bA^1})=0$. By Lemma \ref{interpolation computation}, we have
$$\Fut(\mY, c\mG;-K_{\mY/\bA^1})=c\cdot\Fut(\mY,\mG;-K_{\mY/\bA^1})+(1-c)\cdot\Fut(\mY;-K_{\mY/\bA^1}) \geq 0.$$
Thus we see the following
$$\Fut(\mY; -K_{\mY/\bA^1})\stackrel{\textit{up to a positive multiple}}{=}A_Y(E)-S_Y(E)\geq 0.$$

It remains to prove $(3)\Rightarrow (2)$. Suppose it is not ture, then there exists a sequence of increasing rational numbers $0<a_i<1$ tending to one and a corresponding sequence of lc Calabi-Yau pairs $(Y_i,G_i)$ such that 
\begin{enumerate}
\item $Y_i$ is a $\bQ$-Fano variety of dimension $d$ and $G_i\sim_\bQ -K_{Y_i}$ has its coefficients contained in the set of coefficients of $G$,
\item $A_{Y_i}(F)-S_{Y_i}(F)\geq 0$ for any lc place $F$ of $(Y_i, G_i)$,
\item $(Y_i, a_i G_i)$ is K-unstable for each $i$.
\end{enumerate}
By \cite{LX14}, there exists a special divisor $E_i$ over $Y_i$ which corresponds to a special test configuration of $(Y_i,a_iG_i)$, denoted by $(\mY_i,a_i\mG_i)$, such that 
$$\Fut(\mY_i,a_i\mG_i;-K_{\mY_i})\stackrel{\textit{up to a positive multiple}}{=}A_{Y_i,a_iG_i}(E_i)-S_{Y_i,a_iG_i}(E_i)<0.$$
After subtracting those $a_i$ which are not close to one, we may assume $(\mY_{i,0}, \mG_{i,0})$ is log canonical by \cite{HMX14}. That is to say, the test configuration $(\mY_i,\mG_i)$ degenerates an lc pair $(Y_i,G_i)$ to another lc pair $(\mY_{i,0},\mG_{i,0})$, which implies that $E_i$ is an lc place of $(Y_i,G_i)$ by Lemma \ref{lemma on cy}. This is a contradiction, since $A_{Y_i}(F)-S_{Y_i}(F)\geq 0$ for any lc place $F$ of $(Y_i,G_i)$, however, 
$$A_{Y_i,a_iG_i}(E_i)-S_{Y_i,a_iG_i}(E_i)=(1-a_i)(A_{Y_i}(E_i)-S_{Y_i}(E_i))<0.$$
The proof is finished.
\end{proof}

If $r\geq 1$, we choose $t'_0$ as above; if $0<r<1$, we take $t'_0= 1-r$. From now on, we consider a set $\mF$ of pairs satisfying that $(X,D)\in \mF$ if and only if it is  of the following form: 
\begin{enumerate}
\item $X$ is a $\bQ$-Fano variety of dimension $d$ and volume $v$,
\item $D\sim_\bQ -rK_X$ is a Weil divisor on $X$ such that $(X, \frac{c}{r}D)$ is K-semistable for some rational $c\in [0,1)$.
\end{enumerate}

\begin{theorem}\label{bounded lemma}
Notation as above, the set $\mF$ is contained in a log bounded family.
\end{theorem}

\begin{proof}
Suppose $(X,D)\in \mF$, then $\vol(-K_X)=v$. By Lemma \ref{almost cy}, $(X,\frac{c}{r}D)$ is K-semistable for some rational $c\in [0,1-t_0']$. By valuative criterion, we have
$$\frac{A_X(E)}{S_X(E)} \geq \frac{(1-c)\cdot A_{X, \frac{c}{r}D}(E)}{ S_{X, \frac{c}{r}D}(E) }\geq t_0',$$
where $E$ runs through all prime divisors over $X$.
Thus $\delta(X)\geq t'_0$. By \cite{Jiang20} we know that $X$ lies in a bounded family. As $D\sim_\bQ -rK_X$ is a Weil divisor, the degree of $D$ is also bounded. Thus $\mF$ is contained in a log bounded family.
\end{proof}

Note that the varieties in $\mF$ admit dimension $d$ and volume $v$. We choose a sufficiently divisible $m$ depending only on $(d,v,t'_0)$ such that $-mK_X$ is very ample for every $\bQ$-Fano variety $X$ with $\dim (X)=d$, $\vol(-K_X)=v$, and $\delta(X)\geq t_0'$. Then we consider the following subset of $\mF$, denoted by $\mF_{m,\bP^{N_m}}^{\chi,\tilde{\chi}}$, such that $(X,D)\in \mF_{m,\bP^{N_m}}^{\chi,\tilde{\chi}}$ if and only if the following conditions are satisfied:
\begin{enumerate}
\item $(X,D)$ can be embedded into $\bP^{N_m}\times \bP^{N_m}$\footnote{By this notation we mean that $X$ and $D$ are embedded into $\bP^{N_m}$ respectively.} such that $\mO_X(1):=\mO_{\bP^{N_m}}(1)|_X\sim -mK_X$ and $H^0(\bP^{N_m}, \mO_{\bP^{N_m}}(1))\cong H^0(X,-mK_X)$,
\item $X$ (resp. $D$) admits the Hilbert polynomial $\chi$ (resp. $\tilde{\chi}$)
with respect to $\mO_X(1)$ (resp. $\mO_D(1):=\mO_{\bP^{N_m}}(1)|_D$).
\end{enumerate}
We just use $\mF_m$ to denote $\mF_{m,\bP^{N_m}}^{\chi,\tilde{\chi}}$ for simplicity, and we from now on fix  $P:=(\chi,\tilde{\chi})$. Define the subset $\mF_{m,c}\subset \mF_m$ for a rational $c\in [0,1)$ as follows:$$\mF_{m,c}:=\{(X,D)\in \mF_m| \textit{\rm the pair $(X,\frac{c}{r}D)$ is K-semistable}\}. $$
By the proof of \cite[Theorem 2.21]{XZ20b}, $\mF_m$ admits an algebraic structure, and $\mF_{m,c}\subset \mF_m$ is an open sub-scheme (e.g.  \cite{BLX19,Xu20}), then we write $\mF_m^{\red}$ and $\mF_{m,c}^{\red}$
to be the corresponding reduced schemes. By \cite[Theorem 3.24]{ADL19},  the reduced Artin stacks 
$[\mF_{m,c}^{\red}/\PGL(N_m+1)]$ stabilize for sufficiently divisible $m$, thus we just denote it by $\mM^K_{P,r,c}$. By the works \cite{Jiang20, BLX19, Xu20}, we currently know that $\mM^K_{P,r,c}$ is a reduced Artin stack of finite type representing the moduli pseudo-functor over reduced bases, where $\mM^K_{P,r,c}(S)$ for a reduced scheme $S\in Sch_\bC$ consists of all families $(\mX,\mD)\to S$ satisfying the following conditions:
\begin{enumerate}
\item $(\mX, \mD)\to S$ is a $\bQ$-Gorenstein flat family with each fiber contained in $\mF_{m,c}$,
\item $-K_{\mX/S}$ is a relatively ample $\bQ$-line bundle on $\mX$ and $\mD\sim_{\bQ,S} -rK_{\mX/S}$.
\end{enumerate}

\begin{remark}\label{moduli functor}
We only consider the above moduli pseudo-functor for convenience, but it is also defined for any base (including non-reduced ones) due to the work \cite{Kollar19}. Fix a rational number $c\in [0,1)$, to define the moduli functor for an arbitrary base scheme $S\in Sch_\bC$, one just consider all the families $(\mX,\mD)\to S$ satisfying the following conditions:
\begin{enumerate}
\item $\mX\to S$ is proper and flat,
\item $\mD$ is a K-flat family of divisors on $\mX$ (see \cite{Kollar19}),
\item $-K_{\mX/S}-c\mD$ is $\bQ$-Cartier with each geometric fiber contained in $\mF_{m,c}$.

\end{enumerate}
\end{remark}

By the works \cite{ABHLX20, BHLLX21, XZ20b, LXZ21}, $\mM^K_{P,r,c}$ admits a projective separated good moduli space $M^K_{P,r,c}$ parametrizing K-polystable elements.
We denote $(\mX_m, \mD_m)\subset (\mF_m^{\red}\times \bP^{N_m},\mF_m^{\red} \times \bP^{N_m})$ to be the universal family corresponding to $\mF_m^{\red}$, and denote by 
$$p_1: (\mX_m,\frac{1}{r}\mD_m)=:(\mX, \frac{1}{r}\mD)\to \mF^{\red}_m.$$
Combine Theorem \ref{thm: wall-crossing} and Lemma \ref{almost cy}, we have the following result.

\begin{theorem}\label{wall crossing}
There exist finite rational numbers $0=c_0<c_1<c_2<...<c_k<c_{k+1}=1$ such that 
\begin{enumerate}
\item the schemes $\mF_{m,c}^{\red}$ do not change as $c$ varies in $(c_j,c_{j+1})$ for $0\leq j\leq k$, and we have open immersions
$$\mF^{\red}_{m,c_j-\epsilon}\hookrightarrow \mF^{\red}_{m,c_j}\hookleftarrow \mF^{\red}_{m,c_j+\epsilon}$$
for rational $0<\epsilon\ll1$ and $1\leq j\leq k$. 
\item the Artin stacks $\mM^K_{P,r,c}$  do not change as $c$ varies in $(c_j,c_{j+1})$ for $0\leq j\leq k$, and we have open immersions 
$$\mM^K_{P,r,c_j-\epsilon}\hookrightarrow \mM^K_{P,r,c_j}\hookleftarrow  \mM^K_{P,r,c_j+\epsilon}$$
for rational $0<\epsilon\ll1$ and $1\leq j\leq k$. 
\item  the good moduli spaces $M^K_{P,r,c}$ do not change as $c$ varies in $(c_j,c_{j+1})$ for $0\leq j\leq k$, and we have the following projective morphisms among good moduli spaces
$$M^K_{P,r,c_j-\epsilon}\stackrel{\phi^-_j}{\rightarrow}M^K_{P,r,c_j}\stackrel{\phi^+_j}{\leftarrow}M^K_{P,r,c_j+\epsilon}$$
for rational $0<\epsilon\ll 1$ and $1\leq j\leq k$, which are induced by the following diagram:
\begin{center}
	\begin{tikzcd}[column sep = 2em, row sep = 2em]
	 \mM_{P,r, c_j-\epsilon}^K \arrow[d,"",swap] \arrow[rr,""]&& \mM_{P,r, c_j}^K\arrow[d,"",swap] &&\mM_{P,r, c_j+\epsilon}^K\arrow[d,""]\arrow[ll,""]\\
	 M_{P,r, c_j-\epsilon}^K\arrow[rr,"\phi_j^-"]&& M_{P,r, c_j}^K&&M_{P,r, c_j+\epsilon}^K\arrow[ll,"\phi_j^+", swap].
	 	 	 	\end{tikzcd}
\end{center}
\end{enumerate}
\end{theorem}

\begin{remark}

Consider the universal family 
$$p_1: (\mX_m,\frac{1}{r}\mD_m)=:(\mX, \frac{1}{r}\mD)\to \mF^{\red}_m.$$
Take $0=l_0<l_1<...<l_r<l_{r+1}=1$ and $0=u_0<u_1<...<u_s<u_{s+1}=1$ as in Section \ref{sec: section 3}. Combine Lemma \ref{detailed stability}, we have the following possibilities.
\begin{enumerate}
\item For $1\leq j\leq k$, if $c_j=l_p$ for some $1\leq p\leq r$ but $c_j\neq u_q$ for any $1\leq q\leq s$. By Lemma \ref{detailed stability}, the set $\{(\mX_t,\frac{c_j}{r}\mD_t)| l(t)=c_j\}$ contains no K-polystable points. Thus, $(\mX_t, \frac{c_j}{r}\mD_t)$ is K-polystable if and only if $(\mX_t, \frac{c_j-\epsilon}{r}\mD_t)$ is K-polystable for $0<\epsilon\ll 1$.
This implies that 
$$\phi^-_j: M^K_{P,r,c_j-\epsilon}\to M^K_{P,r,c_j}$$ 
is an isomorphism. %For the morphism $\phi^+_j: M^K_{P,r,c_j+\epsilon}\to M^K_{P,r,c_j}$, the exceptional locus consist of those points corresponding to pairs of the form $(\mX_t, \mD_t)$ such that $(\mX_t, \frac{c_j+\epsilon}{r}\mD_t)$ is K-polystable while $(\mX_t, \frac{c_j}{r}\mD_t)$ is strictly K-semisatble.
\item For $1\leq j\leq k$, if $c_j=u_q$ for some $1\leq q\leq s$ but $c_j\neq l_p$ for any $1\leq p\leq r$, the same analysis as in the previous case tells us that   
$$\phi^+_j: M^K_{P,r,c_j+\epsilon}\to M^K_{P,r,c_j}$$ 
is an isomorphism for $0<\epsilon\ll 1$.
%For the morphism $\phi^-_j: M^K_{P,r,c_j-\epsilon}\to M^K_{P,r,c_j}$, the exceptional locus consist of those points corresponding to pairs of the form $(\mX_t, \mD_t)$ such that $(\mX_t, \frac{c_j-\epsilon}{r}\mD_t)$ is K-polystable while $(\mX_t, \frac{c_j}{r}\mD_t)$ is strictly K-semisatble.
%\item For $1\leq j\leq k$, if $c_j=l_p=u_q$ for some $1\leq p\leq r$ and $1\leq q\leq s$, combine the previous two cases we could characterize the exceptional locus of $\phi_j^+$ and $\phi_j^-$.
\end{enumerate}
\end{remark}

\begin{remark}
By Theorem \ref{thm: wall-crossing} (1), we also have the following open embeddings for rational $0<\epsilon\ll 1$:
$$\mF^{\red}_{m,c_0}\hookleftarrow \mF^{\red}_{m,c_0+\epsilon}\quad \text{and}\quad\mM^K_{P,r,c_0}\hookleftarrow  \mM^K_{P,r,c_0+\epsilon}, $$
which induce the projective morphism
$$\phi_0^+:  M^K_{P,r,c_0+\epsilon}\to M^K_{P,r,c_0}$$
via the following commutative diagram:
\begin{center}
	\begin{tikzcd}[column sep = 2em, row sep = 2em]
	 \mM^K_{P,r,c_0+\epsilon} \arrow[d,"", swap] \arrow[rr,""]&& \mM^K_{P,r,c_0} \arrow[d,"",swap]\\
	 M^K_{P,r,c_0+\epsilon}\arrow[rr,"\phi_0^+"]&& M^K_{P,r,c_0}.
	\end{tikzcd}
\end{center}
By  \cite[Theorem 1.1]{Zhou21a}, if $X$ is a K-polystable Fano variety corresponding to a closed point $[X]\in M^K_{P,r,c_0}$, then the set of points corresponding to pairs of the form $(X, D)$, where $D\sim_\bQ-rK_X$ is a Weil divisor and is GIT-polystable under $\Aut(X)$-action (see \cite{Zhou21a}), is contained in the fiber of $\phi_0^+$ over $[X]$.
\end{remark}

\bibliography{reference.bib}

\begin{bibdiv}
\begin{biblist}

\bib{ABHLX20}{article}{
      author={Alper, Jarod},
      author={Blum, Harold},
      author={Halpern-Leistner, Daniel},
      author={Xu, Chenyang},
       title={Reductivity of the automorphism group of {K}-polystable {F}ano
  varieties},
        date={2020},
        ISSN={0020-9910},
     journal={Invent. Math.},
      volume={222},
      number={3},
       pages={995\ndash 1032},
         url={https://doi.org/10.1007/s00222-020-00987-2},
      review={\MR{4169054}},
}

\bib{ADL19}{article}{
      author={Ascher, Kenneth},
      author={DeVleming, Kristin},
      author={Liu, Yuchen},
       title={Wall crossing for {K}-moduli spaces of plane curves},
        date={2019},
        note={in preparation},
}

\bib{BCHM10}{article}{
      author={Birkar, Caucher},
      author={Cascini, Paolo},
      author={Hacon, Christopher~D.},
      author={McKernan, James},
       title={Existence of minimal models for varieties of log general type},
        date={2010},
        ISSN={0894-0347},
     journal={J. Amer. Math. Soc.},
      volume={23},
      number={2},
       pages={405\ndash 468},
         url={https://doi.org/10.1090/S0894-0347-09-00649-3},
      review={\MR{2601039}},
}

\bib{BHJ17}{article}{
      author={Boucksom, S\'{e}bastien},
      author={Hisamoto, Tomoyuki},
      author={Jonsson, Mattias},
       title={Uniform {K}-stability, {D}uistermaat-{H}eckman measures and
  singularities of pairs},
        date={2017},
        ISSN={0373-0956},
     journal={Ann. Inst. Fourier (Grenoble)},
      volume={67},
      number={2},
       pages={743\ndash 841},
         url={http://aif.cedram.org/item?id=AIF_2017__67_2_743_0},
      review={\MR{3669511}},
}

\bib{BHLLX21}{article}{
      author={Blum, Harold},
      author={Halpern-Leistner, Daniel},
      author={Liu, Yuchen},
      author={Xu, Chenyang},
       title={On properness of {K}-moduli spaces and optimal degenerations of
  {F}ano varieties},
        date={2021},
        ISSN={1022-1824},
     journal={Selecta Math. (N.S.)},
      volume={27},
      number={4},
       pages={Paper No. 73, 39},
         url={https://doi.org/10.1007/s00029-021-00694-7},
      review={\MR{4292783}},
}

\bib{Birkar19}{article}{
      author={Birkar, Caucher},
       title={Anti-pluricanonical systems on {F}ano varieties},
        date={2019},
        ISSN={0003-486X},
     journal={Ann. of Math. (2)},
      volume={190},
      number={2},
       pages={345\ndash 463},
         url={https://doi.org/10.4007/annals.2019.190.2.1},
      review={\MR{3997127}},
}

\bib{BJ20}{article}{
      author={Blum, Harold},
      author={Jonsson, Mattias},
       title={Thresholds, valuations, and {K}-stability},
        date={2020},
        ISSN={0001-8708},
     journal={Adv. Math.},
      volume={365},
       pages={107062, 57},
         url={https://doi.org/10.1016/j.aim.2020.107062},
      review={\MR{4067358}},
}

\bib{BLX19}{article}{
      author={Blum, Harold},
      author={Liu, Yuchen},
      author={Xu, Chenyang},
       title={{Openness of K-semistability for Fano varieties}},
        date={2022},
     journal={Duke Mathematical Journal},
      volume={171},
      number={13},
       pages={2753 \ndash  2797},
         url={https://doi.org/10.1215/00127094-2022-0054},
}

\bib{BLZ19}{misc}{
      author={Blum, Harold},
      author={Liu, Yuchen},
      author={Zhou, Chuyu},
       title={Optimal destabilization of {K}-unstable {F}ano varieties via
  stability thresholds},
   publisher={arXiv},
        date={2019},
         url={https://arxiv.org/abs/1907.05399},
}

\bib{BX19}{article}{
      author={Blum, Harold},
      author={Xu, Chenyang},
       title={Uniqueness of {${K}$}-polystable degenerations of {F}ano
  varieties},
        date={2019},
        ISSN={0003-486X},
     journal={Ann. of Math. (2)},
      volume={190},
      number={2},
       pages={609\ndash 656},
         url={https://doi.org/10.4007/annals.2019.190.2.4},
      review={\MR{3997130}},
}

\bib{CDS15a}{article}{
      author={Chen, Xiuxiong},
      author={Donaldson, Simon},
      author={Sun, Song},
       title={K\"{a}hler-{E}instein metrics on {F}ano manifolds. {I}:
  {A}pproximation of metrics with cone singularities},
        date={2015},
        ISSN={0894-0347},
     journal={J. Amer. Math. Soc.},
      volume={28},
      number={1},
       pages={183\ndash 197},
         url={https://doi.org/10.1090/S0894-0347-2014-00799-2},
      review={\MR{3264766}},
}

\bib{CDS15b}{article}{
      author={Chen, Xiuxiong},
      author={Donaldson, Simon},
      author={Sun, Song},
       title={K\"{a}hler-{E}instein metrics on {F}ano manifolds. {II}: {L}imits
  with cone angle less than {$2\pi$}},
        date={2015},
        ISSN={0894-0347},
     journal={J. Amer. Math. Soc.},
      volume={28},
      number={1},
       pages={199\ndash 234},
         url={https://doi.org/10.1090/S0894-0347-2014-00800-6},
      review={\MR{3264767}},
}

\bib{CDS15c}{article}{
      author={Chen, Xiuxiong},
      author={Donaldson, Simon},
      author={Sun, Song},
       title={K\"{a}hler-{E}instein metrics on {F}ano manifolds. {III}:
  {L}imits as cone angle approaches {$2\pi$} and completion of the main proof},
        date={2015},
        ISSN={0894-0347},
     journal={J. Amer. Math. Soc.},
      volume={28},
      number={1},
       pages={235\ndash 278},
         url={https://doi.org/10.1090/S0894-0347-2014-00801-8},
      review={\MR{3264768}},
}

\bib{CP21}{article}{
      author={Codogni, Giulio},
      author={Patakfalvi, Zsolt},
       title={Positivity of the {CM} line bundle for families of {K}-stable klt
  {F}ano varieties},
        date={2021},
        ISSN={0020-9910},
     journal={Invent. Math.},
      volume={223},
      number={3},
       pages={811\ndash 894},
         url={https://doi.org/10.1007/s00222-020-00999-y},
      review={\MR{4213768}},
}

\bib{Don02}{article}{
      author={Donaldson, S.~K.},
       title={Scalar curvature and stability of toric varieties},
        date={2002},
        ISSN={0022-040X},
     journal={J. Differential Geom.},
      volume={62},
      number={2},
       pages={289\ndash 349},
         url={http://projecteuclid.org/euclid.jdg/1090950195},
      review={\MR{1988506}},
}

\bib{FO18}{article}{
      author={Fujita, Kento},
      author={Odaka, Yuji},
       title={On the {K}-stability of {F}ano varieties and anticanonical
  divisors},
        date={2018},
        ISSN={0040-8735},
     journal={Tohoku Math. J. (2)},
      volume={70},
      number={4},
       pages={511\ndash 521},
         url={https://doi.org/10.2748/tmj/1546570823},
      review={\MR{3896135}},
}

\bib{Fuj19}{article}{
      author={Fujita, Kento},
       title={A valuative criterion for uniform {K}-stability of
  {$\Bbb{Q}$}-{F}ano varieties},
        date={2019},
     journal={J. Reine Angew. Math.},
      volume={751},
       pages={309\ndash 338},
}

\bib{HMX14}{article}{
      author={Hacon, Christopher~D.},
      author={McKernan, James},
      author={Xu, Chenyang},
       title={A{CC} for log canonical thresholds},
        date={2014},
        ISSN={0003-486X},
     journal={Ann. of Math. (2)},
      volume={180},
      number={2},
       pages={523\ndash 571},
         url={https://doi.org/10.4007/annals.2014.180.2.3},
      review={\MR{3224718}},
}

\bib{Jiang20}{article}{
      author={Jiang, Chen},
       title={Boundedness of {$\Bbb Q$}-{F}ano varieties with degrees and
  alpha-invariants bounded from below},
        date={2020},
        ISSN={0012-9593},
     journal={Ann. Sci. \'{E}c. Norm. Sup\'{e}r. (4)},
      volume={53},
      number={5},
       pages={1235\ndash 1248},
         url={https://doi.org/10.24033/asens.244},
      review={\MR{4174851}},
}

\bib{KM98}{book}{
      author={Koll\'{a}r, J\'{a}nos},
      author={Mori, Shigefumi},
       title={Birational geometry of algebraic varieties},
      series={Cambridge Tracts in Mathematics},
   publisher={Cambridge University Press, Cambridge},
        date={1998},
      volume={134},
        ISBN={0-521-63277-3},
         url={https://doi.org/10.1017/CBO9780511662560},
        note={With the collaboration of C. H. Clemens and A. Corti, Translated
  from the 1998 Japanese original},
      review={\MR{1658959}},
}

\bib{Kollar13}{book}{
      author={Koll\'{a}r, J\'{a}nos},
       title={Singularities of the minimal model program},
      series={Cambridge Tracts in Mathematics},
   publisher={Cambridge University Press, Cambridge},
        date={2013},
      volume={200},
        ISBN={978-1-107-03534-8},
         url={https://doi.org/10.1017/CBO9781139547895},
        note={With a collaboration of S\'{a}ndor Kov\'{a}cs},
      review={\MR{3057950}},
}

\bib{Kollar19}{misc}{
      author={Kollár, János},
       title={Families of divisors},
        date={2019},
}

\bib{Li17}{article}{
      author={Li, Chi},
       title={K-semistability is equivariant volume minimization},
        date={2017},
        ISSN={0012-7094},
     journal={Duke Math. J.},
      volume={166},
      number={16},
       pages={3147\ndash 3218},
         url={https://doi.org/10.1215/00127094-2017-0026},
      review={\MR{3715806}},
}

\bib{LWX19}{article}{
      author={Li, Chi},
      author={Wang, Xiaowei},
      author={Xu, Chenyang},
       title={On the proper moduli spaces of smoothable {K}\"{a}hler-{E}instein
  {F}ano varieties},
        date={2019},
        ISSN={0012-7094},
     journal={Duke Math. J.},
      volume={168},
      number={8},
       pages={1387\ndash 1459},
         url={https://doi.org/10.1215/00127094-2018-0069},
      review={\MR{3959862}},
}

\bib{LX14}{article}{
      author={Li, Chi},
      author={Xu, Chenyang},
       title={Special test configuration and {K}-stability of {F}ano
  varieties},
        date={2014},
        ISSN={0003-486X},
     journal={Ann. of Math. (2)},
      volume={180},
      number={1},
       pages={197\ndash 232},
         url={https://doi.org/10.4007/annals.2014.180.1.4},
      review={\MR{3194814}},
}

\bib{LXZ21}{article}{
      author={Liu, Yuchen},
      author={Xu, Chenyang},
      author={Zhuang, Ziquan},
       title={Finite generation for valuations computing stability thresholds
  and applications to {K}-stability},
        date={2022},
        ISSN={0003-486X},
     journal={Ann. of Math. (2)},
      volume={196},
      number={2},
       pages={507\ndash 566},
         url={https://doi.org/10.4007/annals.2022.196.2.2},
      review={\MR{4445441}},
}

\bib{Odaka13b}{article}{
      author={Odaka, Yuji},
       title={A generalization of the {R}oss-{T}homas slope theory},
        date={2013},
        ISSN={0030-6126},
     journal={Osaka J. Math.},
      volume={50},
      number={1},
       pages={171\ndash 185},
         url={http://projecteuclid.org/euclid.ojm/1364390425},
      review={\MR{3080636}},
}

\bib{Oda13}{article}{
      author={Odaka, Yuji},
       title={The {GIT} stability of polarized varieties via discrepancy},
        date={2013},
        ISSN={0003-486X},
     journal={Ann. of Math. (2)},
      volume={177},
      number={2},
       pages={645\ndash 661},
         url={https://doi.org/10.4007/annals.2013.177.2.6},
      review={\MR{3010808}},
}

\bib{Tian15}{article}{
      author={Tian, Gang},
       title={K-stability and {K}\"{a}hler-{E}instein metrics},
        date={2015},
        ISSN={0010-3640},
     journal={Comm. Pure Appl. Math.},
      volume={68},
      number={7},
       pages={1085\ndash 1156},
         url={https://doi.org/10.1002/cpa.21578},
      review={\MR{3352459}},
}

\bib{Tian97}{article}{
      author={Tian, Gang},
       title={K\"{a}hler-{E}instein metrics with positive scalar curvature},
        date={1997},
        ISSN={0020-9910},
     journal={Invent. Math.},
      volume={130},
      number={1},
       pages={1\ndash 37},
         url={https://doi.org/10.1007/s002220050176},
      review={\MR{1471884}},
}

\bib{TW20}{article}{
      author={Tian, Gang},
      author={Wang, Feng},
       title={On the existence of conic {K}\"{a}hler-{E}instein metrics},
        date={2020},
        ISSN={0001-8708},
     journal={Adv. Math.},
      volume={375},
       pages={107413, 42},
         url={https://doi.org/10.1016/j.aim.2020.107413},
      review={\MR{4170229}},
}

\bib{Wang12}{article}{
      author={Wang, Xiaowei},
       title={Height and {GIT} weight},
        date={2012},
        ISSN={1073-2780},
     journal={Math. Res. Lett.},
      volume={19},
      number={4},
       pages={909\ndash 926},
         url={https://doi.org/10.4310/MRL.2012.v19.n4.a14},
      review={\MR{3008424}},
}

\bib{Xu20}{article}{
      author={Xu, Chenyang},
       title={A minimizing valuation is quasi-monomial},
        date={2020},
        ISSN={0003-486X},
     journal={Ann. of Math. (2)},
      volume={191},
      number={3},
       pages={1003\ndash 1030},
         url={https://doi.org/10.4007/annals.2020.191.3.6},
      review={\MR{4088355}},
}

\bib{XZ20b}{article}{
      author={Xu, Chenyang},
      author={Zhuang, Ziquan},
       title={On positivity of the {CM} line bundle on {K}-moduli spaces},
        date={2020},
        ISSN={0003-486X},
     journal={Ann. of Math. (2)},
      volume={192},
      number={3},
       pages={1005\ndash 1068},
         url={https://doi.org/10.4007/annals.2020.192.3.7},
      review={\MR{4172625}},
}

\bib{Zhou21a}{misc}{
      author={Zhou, Chuyu},
       title={Log {K}-stability of {GIT}-stable divisors on {F}ano manifolds},
   publisher={arXiv},
        date={2021},
         url={https://arxiv.org/abs/2102.07458},
}

\bib{ZZ19}{article}{
      author={Zhou, Chuyu},
      author={Zhuang, Ziquan},
       title={Some criteria for uniform {K}-stability},
        date={2021},
        ISSN={1073-2780},
     journal={Math. Res. Lett.},
      volume={28},
      number={5},
       pages={1613\ndash 1632},
      review={\MR{4471722}},
}

\end{biblist}
\end{bibdiv}
\end{document}